\documentclass{compositio}     

\usepackage{amsmath}            
\usepackage[all]{xy}

\begin{document}


\def\M#1{\mathbb#1}     
\def\B#1{\bold#1}       
\def\C#1{\mathcal#1}    
\def\E#1{\scr#1}        

\def\mR{\M{R}}           
\def\mZ{\M{Z}}           
\def\mN{\M{N}}           
\def\mQ{\M{Q}}       
\def\mC{\M{C}}  
\def\mF{\M{F}} 
\def\mP{\M{P}}
\def\mG{\M{G}}



\def\Spec{{\rm Spec}}
\def\PGL{{\rm PGL}}
\def\rg{{\rm rg}}
\def\Hom{{\rm Hom}}
\def\Aut{{\rm Aut}}
 \def\Tr{{\rm Tr}}
 \def\Exp{{\rm Exp}}
 \def\Gal{{\rm Gal}}
 \def\End{{\rm End}}
 \def\det{{{\rm det}}}
 \def\Td{{\rm Td}}
 \def\ch{{\rm ch}}
 \def\che{{\rm ch}_{\rm eq}}
  \def\Spec{{\rm Spec}}
\def\Id{{\rm Id}}
\def\Zar{{\rm Zar}}
\def\Supp{{\rm Supp}}
\def\eq{{\rm eq}}
\def\Ann{{\rm Ann}}
\def\LT{{\rm LT}}
\def\CT{{\rm CT}}
\def\mn{\mu_n}
\def\ac{{\rm ac}}


 \def\pr{^{\prime}}
 \def\prpr{^{\prime\prime}}
 \def\mtr#1{\overline{#1}}
 \def\ra{\rightarrow}
 
 \def\mQ{{\Bbb Q}}
 \def\mR{{\Bbb R}}
 \def\mZ{{\Bbb Z}}
 \def\mC{{\Bbb C}}
 \def\mN{{\Bbb N}}
 \def\CI{{\cal I}}
 \def\CH{{\cal H}}
 \def\CO{{\cal O}}
 \def\CA{{\cal A}}
 \def\CB{{\cal B}}
 \def\CL{{\cal L}}
 \def\CC{{\cal C}}
 \def\CM{{\cal M}}
  \def\CE{{\cal E}}
 
 \def\refeq#1{(\ref{#1})}
 \def\blb{{\big(}}
 \def\brb{{\big)}}
\def\mc{{{\mathfrak c}}}
\def\mcpr{{{\mathfrak c}'}}
\def\mcprpr{{{\mathfrak c}''}}
\def\ul#1{\overline{#1}}
\def\ss{{\rm ss}}
\def\parf{{\rm parf}}
\def\P1{{{\bf P}^1}}
\def\cod{{\rm cod}}
\def\pr{\prime}
\def\prpr{\prime\prime}
\def\ss{\scriptstyle}
\def\OX{{ {\cal O}_X}}
\def\mpartial{{\mtr{\partial}}}
\def\inv{{\rm inv}}
\def\indlim{\underrightarrow{\lim}}
\def\prolim{\underleftarrow{\lim}}
\def\pprolim{'\prolim'}
\def\Pro{{\rm Pro}}
\def\Ind{{\rm Ind}}
\def\Ens{{\rm Ens}}
\def\without{\backslash}
\def\pbdb{{\Pro_b\ D^-_c}}
\def\qc{{\rm qc}}
\def\Com{{\rm Com}}
\def\an{{\rm an}}
\def\gfield{{\rm\bf k}}
\def\red{{\rm red}}
\def\rk{{\rm rk}}
\def\Pic{{\rm Pic}}
\def\gcd{{\rm gcd}}
\def\deg{{\rm deg}}
\def\lb{{\rm lb}}
\def\ell{{\rm ell}}


\def\RHom{{\rm RHom}}
\def\rRHom{{\mathcal RHom}}
\def\rHom{{\mathcal Hom}}
\def\dotimes{{\overline{\otimes}}}
\def\Ext{{\rm Ext}}
\def\rExt{{\mathcal Ext}}
\def\Tor{{\rm Tor}}
\def\rTor{{\mathcal Tor}}
\def\SP{{\mathfrak S}}

\theoremstyle{plain}
 \newtheorem{theor}{Theorem}[section]
 \newtheorem{prop}[theor]{Proposition}
 \newtheorem{cor}[theor]{Corollary}
 \newtheorem{lemma}[theor]{Lemma}
 \newtheorem{slemma}[theor]{sublemma}
 
 \theoremstyle{definition}
 \newtheorem{defin}[theor]{Definition}
 \newtheorem{conj}[theor]{Conjecture}
 
\theoremstyle{remark} 
\newtheorem{remark}{Remark}

 
\title
  [On the determinant bundles of abelian schemes]         
  {On the determinant bundles of abelian schemes} 
\author{Vincent Maillot}
\email{vmaillot@math.jussieu.fr}  
\address{Institut de Math\'ematiques de Jussieu,
 Universit\'e Paris 7 Denis Diderot, C.N.R.S.,
 Case Postale 7012,
 2 place Jussieu,
 F-75251 Paris Cedex 05, France}
\author{Damian R\"ossler}
\email{dcr@math.jussieu.fr}
\address{Institut de Math\'ematiques de Jussieu,
 Universit\'e Paris 7 Denis Diderot, C.N.R.S.,
 Case Postale 7012,
 2 place Jussieu,
 F-75251 Paris Cedex 05, France}
\shortauthors{Maillot and R\"ossler}
\begin{abstract}
 Let $\pi:\CA\ra S$ be an abelian scheme over a scheme $S$  which is 
 quasi-projective over an affine noetherian scheme 
  and let $\CL$ be a symmetric, 
 rigidified, relatively ample line bundle on $\CA$. We show that there is an isomorphism 
 $$
  \det(\pi_*\CL)^{\otimes 24}\simeq\big(\pi_*\omega_{\CA}^{\vee}\big)^{\otimes 
 12d}
 $$
 of line bundles on $S$, where $d$ is the rank of the (locally free) sheaf $\pi_*\CL$. 
 We also show that the numbers $24$ and $12d$ are sharp in the 
 following sense: if $N>1$ is a common divisor of $12$ and $24$, then 
 there are data as above such that 
 $$
  \det(\pi_*\CL)^{\otimes (24/N)}\not\simeq\big(\pi_*\omega_{\CA}^{\vee}\big)^{\otimes 
 (12d/N)}.
 $$
 \end{abstract}

\maketitle 

 
\section{Introduction}
 
 Let $\pi:\CA\ra S$ be a abelian scheme, where $S$ is a scheme which is 
 quasi-projective over an affine noetherian scheme. We denote as usual by $\omega_{\CA}$ the 
 determinant of the sheaf of differentials of $\pi$. 
  Let $\CL$ be a line bundle on $\CA$. Let $\epsilon:S\ra\CA$ be the 
  zero-section and suppose that the line bundle 
  $\epsilon^*\CL$ is the trivial line bundle. Suppose furthermore that there is an isomorphism 
  $[-1]^*\CL\simeq\CL$ and that $\CL$ is ample relatively to $\pi$. In this situation, 
 Chai and Faltings prove the following result (see \cite[Th. 5.1, p. 25]{Chai-Faltings}):
 \begin{theor}[(Chai-Faltings)]
 There is an isomorphism
 $\det(\pi_*\CL)^{\otimes 8d^3}\simeq\big(\pi_*\omega_{\CA}^{\vee}\big)^{\otimes 
 4d^4}$
 of line bundles on $S$.
 \label{CF}
 \end{theor}
  Here $d$ is the rank of the (locally free) sheaf $\pi_*\CL$. This is 
 a refinement of a special case of the "formule cl\'e" considered by 
 Moret-Bailly in his monograph \cite{Moret-Bailly}.
 
 In \cite[p. 27]{Chai-Faltings}, Chai and Faltings state that it is nevertheless likely that the factor $d^3$ can be cancelled on both 
sides of the above isomorphism or in other words that it is likely that there 
is an isomorphism 
\begin{equation}
 \det(\pi_*\CL)^{\otimes 8}\simeq\big(\pi_*\omega_{\CA}^{\vee}\big)^{\otimes 
 4d}.
 \label{strong}
 \end{equation}
 
 Let us introduce the line bundle
 $$
\Delta(\CL):=\det(\pi_*\CL)^{\otimes 2}\otimes\pi_*\omega_{\CA}^{\otimes d}.
$$
The existence of the isomorphism \refeq{strong} is the statement 
that $\Delta(\CL)^{\otimes 4}$ is trivial.

The aim of this text is to present the proof of the following statements about 
 $\Delta(\CL)$: 
 \begin{theor}
\begin{itemize}
\item[\rm (a)] There is an isomorphism 
\mbox{$\Delta(\CL)^{\otimes 12}\simeq\CO_S$}.
\item[\rm (b)] For every $g\geq 1$, 
there exist data $\pi:\CA\ra S$ and $\CL$ as above such that 
$\dim(\CA/S)=g$ and such that 
$\Delta(\CL)$ 
is of order $12$ in the Picard group of $S$. 
\end{itemize}
\label{MR}
\end{theor}
The following corollary follows immediately from Theorem \ref{MR} (a) and Theorem 
\ref{CF}.
\begin{cor}
 If $(3,d)=1$ then $\Delta(\CL)^{\otimes 4}$ is trivial.
\end{cor}
Notice that Theorem \ref{MR} (b) in particular implies that the exponent $4$ surmised 
by Chai and Faltings is not the right one (it has to be replaced by the exponent $12$). 
The corollary says that the exponent $4$ is nevertheless the right one when 
$(3,d)=1$. 

The  fact that 
$\Delta(\CL)$ is a torsion line bundle 
is a consequence of the Grothendieck-Riemann-Roch theorem. 
This was shown by Moret-Bailly and Szpiro in the Appendix 2 to 
Moret-Bailly's monograph \cite{Moret-Bailly} and also by Chai in his thesis 
(see \cite[Chap. V, par. 3, th. 3.1, p. 209]{Chai}). The link between the 
Grothendieck-Riemann-Roch theorem and the fact that 
 $\Delta(\CL)$ is a torsion line bundle was already known to Mumford 
 in the early sixties (private communication between Chai and the authors). 
If $S$ is a smooth quasi-projective scheme over $\mC$, then 
\ref{MR} (a) is contained in a theorem of Kouvidakis (see \cite[Th. A]{Kouvidakis}). The method 
of proof of the theorem of Kouvidakis is analytic and is based 
on the study of the transformation formulae of theta functions. 
It extends earlier work by Moret-Bailly (see \cite{Moret-Bailly-2}), who considered the case where 
$d=1$. The result of Kouvidakis was extended by Polishchuk 
to more general bases $S$ in \cite{Polishchuk}. Polishchuk's proof is a refinement 
of Chai and Faltings proof of Theorem \ref{CF}; this last proof is not based 
on the Riemann-Roch theorem. The Theorem 0.1 in 
\cite{Polishchuk} shows in particular that there exists a constant $N(g)$, 
which depends only on the relative dimension $g$ of $\CA$ over $S$, 
such that $\Delta(\CL)^{\otimes N(g)}$ is trivial. The 
Theorems 0.1, 0.2, 0.3 of \cite{Polishchuk} give various bounds for $N(g)$, which 
depend on $d$, $g$ and on the residue characteristics of $S$. 
In this context, the content of Theorem \ref{MR} is that $N(g)=12$ 
is a possible choice and that for each $g\geqslant 1$, it is the best possible choice.

A key input in Polishchuk's refinement of the proof of Chai and Faltings is a formula 
describing the behaviour of $\Delta(\CL)$ under isogenies of 
abelian schemes (\cite[Th. 1.1]{Polishchuk}; see also the end of section 2), which generalises an earlier formula 
by Moret-Bailly (see \cite[VIII, 1.1.3, p. 188]{Moret-Bailly}), who considered the case $d=1$. 
The proof of Theorem \ref{MR} (a) presented in this paper combines 
Polishchuk's isogeny formula and a refinement 
of the Grothendieck-Riemann-Roch theorem, called the Adams-Riemann-Roch 
theorem (see section 2). In spirit, it is close to Mumford's original approach. 
 Our method can also be related to Moret-Bailly's proof 
of the "formule cl\'e" in positive characteristic; see the first remark at the end of the text. 
Our proof of Theorem \ref{MR} (b) is based on a lemma of Polishchuk and on 
two constructions of Mumford. 

The plan of the article is as follows. The second section contains some
preliminaries to the proof; these preliminaries are the Adams-Riemann-Roch theorem
and the two results of Polishchuk mentionned in the last paragraph.
The proof itself is contained in the third section. 

{\bf Notation and conventions.} 
Suppose that $\CM$ is a line bundle on a group scheme $\CC$
 over a base $B$, with zero-section 
$\epsilon:B\ra\CC$. We shall say that 
$\CM$ is rigidified if $\epsilon^*\CM$ is the trivial line bundle. We 
shall say that $\CM$ is symmetric, if $[-1]^*\CM\simeq\CM$. 
Suppose that $x$ is an element of an abelian group $G$ and that $k$ is 
a positive integer; we shall say that $x$ is $k^\infty$-torsion element of $G$ if 
there exists an integer $n\geq 0$ such that $k^n\cdot x=0$ in $G$. 
If $G$ is a group or a group functor and $k$ is a strictly positive integer, we shall 
write $[k]$ for the map $G\ra G$ such that $[k](x)=x+\dots+x$ ($k$-times) 
for every $x\in G$.

\section{Preliminaries}

In this section, "scheme" will be short for "noetherian scheme". 

\subsection{The Adams-Riemann-Roch theorem}

 In this section, we first describe the special case of the Adams-Riemann-Roch 
 theorem that we shall need. We then go on to describe Polishchuk's 
 isogeny formula. 
 
 If $Y$ is a scheme, we shall write as usual $K_0(Y)$ for the 
 Grothendieck group of coherent locally free sheaves. The tensor 
 product of locally free sheaves descends to a bilinear pairing on $K_0(Y)$, which makes 
 it into a commutative ring. If $f:X\ra Y$ is a morphism of schemes, 
 the pull-back of $\CO_Y$-modules induces a ring morphism 
 $f^*:K_0(Y)\ra K_0(X)$. As a ring,  
 $K_0(Y)$ is endowed with a family $(\psi^k)_{k\in\mN^*}$ of 
 (ring) endomorphisms, called the Adams operations. They have 
 the property that  $\psi^k(\CM)=\CM^{\otimes k}$ in $K_0(Y)$ 
 for every line bundle $\CM$ on $Y$. Furthermore, if 
 $f:X\ra Y$ is a scheme morphism as before, then 
 $f^*\circ\psi^k=\psi^k\circ f^*$. The Adams operations are uniquely 
 determined by these two last properties and by the fact 
 that they are ring endomorphisms.
 
 We shall also need Bott's "cannibalistic" classes. 
 We shall denote thus a family of operations $(\theta^k)_{k\in\mN^*}$, each of which 
 associates elements of 
 $K_0(Y)$ to coherent locally free sheaves on $Y$. They have the following three 
 properties, which determine them uniquely.  For every line bundle 
 $\CM$ on $Y$, we have 
 $$
 \theta^k(\CM)=1+\CM+\CM^{\otimes 2}+\dots+\CM^{\otimes (k-1)}
 $$
 in $K_0(Y)$. If 
 $$
 0\ra E'\ra E\ra E''\ra 0
 $$
 is an exact sequence of coherent locally free sheaves on $Y$, then 
 $\theta^k(E')\theta^k(E'')=\theta^k(E)$. And finally, if 
 $f:X\ra Y$ is a morphism of schemes, then $f^*(\theta^k(E))=
 \theta^k(f^*E)$, for every coherent locally free sheaf on $Y$.
 About the operations $\theta^k$, the following lemma holds. 
 Suppose for the time of the lemma that $Y$ is quasi-projective 
 over an affine scheme.
 \begin{lemma}
 For any coherent locally free sheaf $E$ on $Y$, the 
 element $\theta^k(E)$ is invertible in the 
 ring $K_0(Y)[{1\over k}]$.
 \end{lemma}
 \begin{proof}
 See \cite[Par. 4, Prop. 4.2]{Roessler} (for lack of a standard reference).
 \end{proof}
 Let $f:X\ra Y$ be a flat and projective morphism of 
 schemes.
  We may consider the Grothendieck 
 group $K_0^\ac(X)$ of $f$-acyclic coherent locally free sheaves on $X$, 
 i.e. coherent locally free sheaves $E$ such that $R^i f_*E=0$ for every $i>0$. There is a unique 
 morphism of groups $f_*:K_0^\ac(X)\ra K_0(Y)$ such that 
 $f_*(E)=R^0 f_*E$ for every coherent locally free sheaf on $X$. 
 A theorem of Quillen (see \cite[Par. 4, Th. 3, p. 108]{Quillen}) 
 now implies that the natural 
 map $K_0^\ac(X)\ra K_0(X)$ is an isomorphism. 
 Hence we obtain a morphism $f_*:K_0(X)\ra K_0(Y)$. 
 This morphism satisfies the projection formula: 
 for all $y\in K_0(Y)$ and all $x\in K_0(X)$, the identity 
 $f_*(f^*(y)\otimes x)=y\otimes f_*(x)$ holds. 
 
 Let us now consider a smooth and projective morphism 
 of schemes $f:X\ra Y$, where $Y$ is quasi-projective over an affine scheme. Let 
 $\Omega$ be the sheaf of differentials associated to $f$; it is a locally
 free sheaf on $X$. 
 In this situation, the Adams-Riemann-Roch theorem is the following statement: 
 \begin{theor}[(Grothendieck et al.)]
 For any $x\in K_0(X)[{1\over k}]$, the equality 
 \begin{equation*}
 \psi^k(f_*(x))=f_*(\theta^k(\Omega)^{-1}\psi^k(x))
 \end{equation*}
 \label{ARR}
 holds in $K_0(Y)[{1\over k}]$.
 \end{theor}
 For a proof of the Adams-Riemann-Roch theorem, see \cite[V, par. 7, Th. 7.6, p. 149]{Fulton-Lang}.

\subsection{Some results of Polishchuk on $\Delta(\CL)$}

Let $T$ be any base scheme and let $\kappa:\CB\ra T$ 
be an abelian scheme. Let $\alpha:\CB\ra\CB$ be a finite 
and flat $T$-homomorphism of group schemes. Let 
$\CM$ be a symmetric and rigidified line bundle on $\CB$, which is ample 
relatively to $\kappa$. Suppose that the sheaf $\kappa_*\CM$ has strictly positive rank. 
The following special case of Polishchuk's isogeny formula \cite[Th. 1.1]{Polishchuk} plays 
a crucial role in the proof of Theorem \ref{MR} a):
\begin{theor}[(Polishchuk)]
\begin{itemize}
\item[\rm (a)]
Let $n:=(12,\deg(\alpha))$. 
There is an isomorphism
\begin{equation*}
 \det(\kappa_*(\alpha^*\CL))^{\otimes 2n}\simeq 
 \det(\kappa_*(\CL))^{\otimes(2n\cdot\deg(\alpha))}.
 \end{equation*}
 \item[\rm (b)]
 Let $m:=(3,\deg(\alpha))$. 
 Suppose that $\deg(\alpha)$ is odd and that $\rk(\kappa_*\CL)$ is even. There is then an isomorphism
 \begin{equation*}
 \det(\kappa_*(\alpha^*\CL))^{\otimes m}\simeq 
 \det(\kappa_*(\CL))^{\otimes(m\cdot\deg(\alpha))}.
 \end{equation*}
 \end{itemize}
  \label{polish}
\end{theor}

The two following lemmata are needed in the proof of Theorem \ref{MR} (b). 
\begin{lemma}[(Polishchuk)]
Suppose that $\dim(\CB/T)=1$ and 
that $\CM=\CO(O_\CB)\otimes\omega_{\CB}$. Then for every $r\geq 1$, there is 
an isomorphism 
$$
\Delta(\CM^{\otimes r})\simeq\omega_{\CB}^{\otimes (r^2+2)}.
$$
\label{polish2}
\end{lemma}
Here $O_\CB$ is the image of the unit section of $\CB/T$. 
Notice that the image of the unit section is a Cartier divisor and that its normal 
bundle is isomorphic to the restriction of $\omega_{\CB}$ via the unit 
section (this is a consequence of the fact that $\kappa$ is smooth; 
see for instance \cite[IV, par. 3, Lemma 3.8]{Fulton-Lang}). This 
implies that the pull-back of $\CM$ via the unit section is the trivial 
line bundle (see \cite[IV, par. 3, Prop. 3.2 (b)]{Fulton-Lang}).

For the proof of the Lemma \ref{polish2}, see \cite[Prop. 5.1]{Polishchuk}. 

Let now $\kappa':\CB'\ra T$ 
be an abelian scheme and let $\CM'$ be a symmetric and rigidified line bundle on $\CB'$, which is ample relatively to $\kappa'$. Let $p$ (resp. $p'$) be the natural projection 
$\CB\times_B\CB'\ra \CB$ (resp. $\CB\times_B\CB'\ra \CB'$). Let 
$m$ (resp. $m'$) be the rank of $\kappa_*\CM$ (resp. $\kappa'_*{\CM'}$). 

\begin{lemma}
There is an isomorphism
$$
\Delta(p^*\CM\otimes {p'}^{*}\CM')\simeq 
\Delta(\CM)^{\otimes m'}\otimes \Delta(\CM')^{\otimes m}
$$
\label{kun}
\end{lemma}
\begin{proof}
Left to the reader (use the K\"unneth formula).
\end{proof}

 \section{The proof}
 
 \subsection{The isomorphism $\Delta(\CL)^{\otimes 12}\simeq\CO_S$}
 
 In this subsection, we shall prove assertion (a) in Theorem \ref{MR}.
 We shall now apply the Adams-Riemann-Roch theorem \ref{ARR} to abelian schemes. 
 We work in the situation of the introduction.   
We may also suppose without 
 restriction of generality that $d\geq 1$. 
  Let $k\geq 2$.  Let $g$ be the relative dimension of $\CA$ over $S$. 
  Recall that the theorem of the cube (see \cite[Par. 5.5, p. 29]{Moret-Bailly}) implies that 
  $[k]^*\CL\simeq\CL^{\otimes k^2}$. Write $\Omega$ for the sheaf of differentials of 
  $\pi$. 
 We compute in $K_0(S)[{1\over k}]$:
 \begin{eqnarray*}
 \psi^{k^2}(\pi_*\CL)&=&\pi_*(\theta^{k^2}(\Omega)^{-1}\psi^{k^2}(\CL))=R\pi_*(\theta^{k^2}(\Omega)^{-1}\CL^{\otimes k^2})\\
 &=&\pi_*(\theta^{k^2}(\Omega)^{-1} [k]^*(\CL))=
 \pi_*([k]^*\CL)\theta^k(\pi_*\Omega)^{-1}
 \end{eqnarray*}
 where we have used the theorem of the cube, 
 the projection formula and the fact that $\pi^*\pi_*\Omega=\Omega$. In other words, 
 we have the identity 
  \begin{eqnarray}
 \theta^{k^2}(\pi_*\Omega)\psi^{k^2}(\pi_*\CL)=
 \pi_*([k]^*\CL)
 \label{mainid}
 \end{eqnarray}
 in $K_0(S)[{1\over k}]$. 
 Let us now introduce the (truncated) Chern character 
 $$
 \ch:K_0(S)[{1\over k}]\ra\mZ[{1\over k}]\oplus\Pic(S)[{1\over k}],
 $$
 which is defined by the formula 
 $$
 \ch(s/k^t):={\rm rank}(s)/k^t\oplus\det(s)^{1/k^t}
 $$
 for every $s\in K_0(S)$ and $t\in\mN$. 
 Let us introduce the pairing
 $$
 (r/k^t,m/k^l)\bullet((r')^{1/k^{t'}},(m')^{1/k^{l'}}):=r\cdot r'/k^{t+t'}\ \oplus\ (m')^{r/k^{t+l'}}\otimes m^{r'/k^{t'+l}}
 $$
 in the group $\mZ[{1\over k}]\oplus\Pic(S)[{1\over k}]$. The pairing $\bullet$
 makes this group into a commutative ring. The properties of the determinant 
 show that the Chern character is a ring morphism. 
 We now apply the Chern character to the identity \refeq{mainid}. As we shall 
 compute in the ring \mbox{$\mZ[{1\over k}]\oplus\Pic(S)[{1\over k}]$}, 
 we switch from multiplicative notation ("$\otimes$") to additive notation ("$+$") in 
 the group $\Pic(S)$. 
 For the purposes of computation, we 
 may suppose without loss of generality that $\pi_*\Omega=
 \omega_1+\dots\omega_g$ in $K_0(S)$, where 
 $\omega_1,\dots,\omega_g$ are line bundles. We compute
 \begin{eqnarray*}
  \ch(\theta^{k^2}(\pi_*\Omega))&=&
  \big(k^2+{k^2(k^2-1)\over 2}\det(\omega_1)\big)\bullet\dots\bullet
  \big(k^2+{k^2(k^2-1)\over 2}\det(\omega_g)\big)\\
  &=&
  k^{2g}+{k^2(k^2-1)k^{2g-2}\over 2}\det(\pi_*\Omega)
  \end{eqnarray*}
  and 
 \begin{eqnarray*}
 \ch(\theta^{k^2}(\pi_*\Omega))\ch(\psi^{k^2}(\pi_*\CL))&=&
 \big(k^{2g}+{k^2(k^2-1)k^{2g-2}\over 2}\det(\pi_*\Omega)\big)\bullet\big(d+
 k^{2}\det(\pi_*\CL)\big)\\
 &=&k^{2g}d+k^{2g+2}\det(\pi_*\CL)+{dk^2(k^2-1)k^{2g-2}\over 2}\det(\pi_*\Omega).
 \end{eqnarray*}
On the other hand, we have  
 $$
 \ch(\pi_*([k]^*\CL))=dk^{2g}+\det(\pi_*[k]^*\CL).
 $$
 Here we have used the fact that the degree of the isogeny given by 
 multiplication by $k$ on $\CA$ is $k^{2g}$ and 
 the fact that the rank of $\pi_*[k]^*\CL$ is $dk^{2g}$ 
 (see \cite[III, par. 12, Th. 2, p. 121]{Mumford}). 
 Thus, \refeq{mainid} leads to the equality 
 \begin{eqnarray*}
k^{2g}d+k^{2g+2}\det(\pi_*\CL)+{dk^2(k^2-1)k^{2g-2}\over 2}\det(\pi_*\Omega)=
dk^{2g}+\det(\pi_*[k]^*\CL)
\end{eqnarray*}
in $\mZ[{1\over k}]\oplus\Pic(S)[{1\over k}]$. Multiplying by 
$k^{-2g}$ and specializing to $\Pic(S)[{1\over k}]$, we get
 \begin{eqnarray*}
k^{2}\det(\pi_*\CL)+{d(k^2-1)\over 2}\det(\pi_*\Omega)=
k^{-2g}\det(\pi_*[k]^*\CL)
\end{eqnarray*}
in $\Pic(S)[{1\over k}]$. Now Theorem \ref{polish} (a) shows that
\begin{equation}
2\cdot k^{-2g}\det(\pi_*[k]^*\CL)=2\cdot\det(\pi_*\CL)
\label{pa}
\end{equation}
in $\Pic(S)[{1\over k}]$. We deduce from the last two equalities that
\begin{eqnarray}
{(k^2-1)}\cdot\big(2\cdot\det(\pi_*\CL)+d\cdot\det(\pi_*\Omega)\big)=0.
\label{finres}
\end{eqnarray}
in $\Pic(S)[{1\over k}]$. In other words, $\Delta(\CL)^{\otimes(k^2-1)}$ is a $k^\infty$-torsion line bundle. 
If we specialise to $k=2$, we see that 
$\Delta(\CL)^{\otimes 3}$ is a $2^\infty$-torsion line bundle. 
If we specialise to $k=3$, we see that 
$\Delta(\CL)^{\otimes 8}$ is a $3^\infty$-torsion line bundle. Hence 
$\Delta(\CL)^{\otimes 24}$ is a  trivial line bundle. 

Suppose now that $d$ is odd. Theorem \ref{CF} says that 
$\Delta(\CL)^{\otimes4d^3}$ is a trivial line bundle. Hence 
$
\Delta(\CL)^{\otimes(24,4d^3)}
$
is a trivial line bundle. Since $(24,4d^3)$ divides $12$, this implies 
that $\Delta(\CL)^{\otimes 12}$ is a  trivial line bundle. 

Suppose now that $d$ is even. Theorem \ref{polish} (b) then shows that  
the equality
\begin{equation*}
k^{-2g}\det(\pi_*[k]^*\CL)=\det(\pi_*\CL)
\end{equation*}
holds in $\Pic(S)[{1\over k}]$ (this equality refines \refeq{pa}). Proceeding 
as we did after the equality \refeq{finres}, we obtain the equality 
\begin{eqnarray*}
{k^2-1\over 2}\big(2\cdot\det(\pi_*\CL)+d\cdot\det(\pi_*\Omega)\big)=0.
\end{eqnarray*}
in $\Pic(S)[{1\over k}]$. In other words, $\Delta(\CL)^{\otimes(k^2-1)}$ is a $k^\infty$-torsion line bundle if $k$ is even and $\Delta(\CL)^{\otimes(k^2-1)/2}$ is 
a $k^\infty$-torsion line bundle if $k$ is odd. If we specialise to $k=3$, we see that 
$\Delta(\CL)^{\otimes 4}$ is a $3^\infty$-torsion line bundle. We saw above 
that $\Delta(\CL)^{\otimes 3}$ is a $2^\infty$-torsion line bundle and so we obtain  
again that $\Delta(\CL)^{\otimes 12}$ is a  trivial line bundle.

This concludes the proof of the assertion (a) of Theorem \ref{MR}.

\begin{remark}
Suppose that $S$ is a scheme over $\mF_p$, for some prime 
number $p$ and that $d=1$. Moret-Bailly then proves that the line bundle 
$\Delta(\CL)^{\otimes(p^2-1)p^{2g+2}}$ is trivial (see \cite[chap. VIII, par. 2, Th. 2.1, p. 193]{Moret-Bailly}). The equality 
\refeq{finres} for $k=p$ is a variant of this. 
Notice furthermore that for any vector bundle 
$E$ on $S$, we have $\psi^p(E)=F_S^*(E)$, where 
$F_S$ is the absolute Frobenius endomorphism of $S$. 
Moret-Bailly's proof is based on the study of the behaviour of $\Delta(\CL)$ under 
base-change by $F_S$ and on the case $d=1$ of the isogeny formula. 
In this sense, our proof of \refeq{mainid} over a general base can be considered as an 
extension of Moret-Bailly's proof of \refeq{mainid} in positive characteristic.
\end{remark}


\subsection{Sharpness}

In this subsection, we shall prove the assertion (b) in 
Theorem \ref{MR}. We fix an affine noetherian base scheme $B$. 
All schemes and morphisms of schemes in this subsection will be 
relative to this base scheme. Furthermore, all schemes will be locally 
noetherian. 
We first recall a result of Mumford.
Consider the following set of data:
\begin{itemize}
\item $\delta,g\in\mN^*$;
\item $T$ a scheme;
\item $\kappa:\CB\ra T$ a projective abelian scheme of relative dimension $g$;
\item $\lambda:\CB\ra \CB^\vee$ a polarisation over $T$ of degree $\delta^2$;
\item a linear rigidification 
$\mP(\kappa_*(L^\Delta(\lambda)^{\otimes 3}))\simeq 
\mP^{6^g\cdot \delta-1}_T$.
\end{itemize}
Here $L^\Delta(\lambda)$ is the 
pull-back of the Poincar\'e line bundle on  $\CB\times_T\CB^\vee$  via 
the map \makebox{$\Id\times_T\lambda:\CB \to\CB \times_T B^\vee$.} 
We shall call the scheme $T$ the ground scheme of the set of data.
If we are given two sets of data as above, there is an obvious notion of 
isomorphism between them. 
If we are given two sets of data with the same ground scheme $T$, an isomorphism 
between the two sets of data will be called a \makebox{$T$-isomorphism} 
if it restricts to the identity on $T$. 
For each scheme $T$, we shall write ${\cal H}_{g,\delta}(T)$ for the set of 
\makebox{$T$-isomorphism} classes of sets of data whose ground scheme is $T$. If 
$T'\ra T$ is a morphism of schemes, the obvious base-change of sets of data from 
$T$ to $T'$ induces a map 
${\cal H}_{g,\delta}(T)\ra{\cal H}_{g,\delta}(T')$. One thus obtains a contravariant 
functor from the category of (locally noetherian) schemes to the category 
of sets. For more details, see \cite[chap. 7, par. 2]{MFK}.
\begin{theor}[(Mumford)] The 
functor ${\cal H}_{g,\delta}$ is representable by a quasi-projective scheme over $B$.
\end{theor}
For the proof, see \cite[Prop. 7.3, chap. 7, par. 2]{MFK}. 
We shall refer to the scheme representing ${\cal H}_{g,\delta}$ as 
${H}_{g,\delta}$. 
 
Let now $\kappa_{1,1}:\CB_{1,1}\ra {H}_{1,1}$ be the universal 
abelian scheme over ${H}_{1,1}$.  Let $$\CL_{\CB_{1,1}}:=\CO(O_{\CB_{1,1}})^{\otimes 3}\otimes\omega^{\otimes 3}_{\CB_{1,1}/H_{1,1}}.$$
 Here again $O_{\CB_{1,1}}$ is the image of the unit section of 
$\CB_{1,1}\ra {H}_{1,1}$. 
\begin{prop}
If $B=\Spec\ \mC$ then 
the line bundle $\Delta(\CL_{\CB_{1,1}})$ is of order $12$ 
in $\Pic({H}_{1,1})$. 
\label{exactor}
\end{prop}
\begin{proof}
Notice that there is a natural action of the group scheme $\PGL_{6}$ 
on $H_{1,1}$, defined as follows. Consider a set of data $\mathfrak D$ of the type described 
at the beginning of the subsection; 
let $a\in \PGL_{6}(T)$; to $a$ corresponds by construction an 
automorphism $A$ of $\mP^{5}_T$; we let $a$ send $\mathfrak D$ on 
the set of data $\mathfrak D$ with its linear rigidification composed with $A$. This defines 
an action of the group functor $\PGL_{6}$ on the functor 
${\cal H}_{1,1}$ and hence an action of $\PGL_{6}$ 
on $H_{1,1}$. Since $H_{1,1}$ is  a fine moduli-space, this 
$\PGL_6$-equivariant structure canonically lifts to a $\PGL_6$-equivariant structure 
on $\CB_{1,1}$, such that the morphism $\kappa_{1,1}$ is 
$\PGL_6$-equivariant. 

Let now $k\geqslant 1$ be the order of $\Delta(\CL_{\CB_{1,1}})$ in $\Pic(H_{1,1})$ 
(which is finite by Theorem \ref{MR} (a)). 
Notice that Lemma \ref{polish2} shows that there is an isomorphism 
$$\Delta(\CL_{\CB_{1,1}})\simeq\kappa_{1,1*}\omega_{\CB}^{11}.$$  We thus see that there is an isomorphism 
\begin{equation}
\kappa_{1,1*}\omega_{\CB_{1,1}}^{\otimes 11\cdot k}\simeq\CO.
\label{ameq}
\end{equation}
Fix such an isomorphism. This is tantamount  to giving 
a trivialising section $s$ of $\kappa_{1,1*}\omega_{\CB_{1,1}}^{\otimes 11\cdot k}$. 
Notice now that the reduced 
closed subscheme $H_{1,1,\red}$ underlying $H_{1,1}$ carries a $\PGL_6$-action 
such that the closed immersion $H_{1,1,\red}\hookrightarrow H_{1,1}$ is 
equivariant (this follows from the definition of a group-scheme action, from 
\cite[I, par. 5, cor. 5.1.8]{EGA} and 
from the fact that $\PGL_6\times H_{1,1}$ is reduced, since $\PGL_6$ is smooth over $\mC$). 
Furthermore, there are no non-trivial characters 
$\PGL_{6}\ra \mG_m$. Thus the restriction of the 
section $s$ to $H_{1,1,\red}$ is 
$\PGL_6$-invariant (for this, see \cite[Prop. 1.4, chap. I, par. 3, p. 33]{MFK}). 

Now consider 
an elliptic curve $\kappa:E\ra\Spec\ \mC$, which has complex multiplication 
by $\mZ[j]$, where $j=-1/2+i\sqrt{3}/2$ is a primitive $3$rd root of unity. 
The element $j$ then acts on $\kappa_*\omega_{E}$ by multiplication by 
either ${j}$ or $\mtr{j}$ (this can be seen by considering the complex uniformisation of $E(\mC)$). 
Choose an arbitrary rigidification of $\mP(\kappa_*(\CO(O_E)^{\otimes 3}))$. 
The elliptic curve $E$ together with its rigidification defines an 
element $P$ of $H_{1,1,\red}(\mC)$, since $\Spec\ \mC$ is reduced.
 Since the section $s$ is  $\PGL_6$-invariant on $H_{1,1,\red}$, 
the element $s(P)\in\kappa_{*}\omega_{E}^{\otimes 11\cdot k}$ must 
satisfy the equation $s(P)={j}^{11\cdot k}\cdot s(P)$ or the equation 
$s(P)=\mtr{j}^{11\cdot k}\cdot s(P)$ and hence 
$3|k$. A similar reasoning with an elliptic curve with 
complex multiplication by the Gaussian integers $\mZ[i]$ shows that 
$4|k$. Hence $12$ divides $k$; on the other hand Theorem \ref{MR} (a) shows that 
$k$ divides $12$. Hence $k=12$ and this concludes the proof.
\end{proof}

\begin{remark}
The idea to use rigidifications in the context of the Theorem \ref{CF} 
is due to Chai and Faltings; see \cite[proof of th. 5.1]{Chai-Faltings}. 
The idea to use elliptic curves with complex multiplication  to compute 
orders in Picard groups is due to Mumford; see \cite[par. 6]{Mumford2}.
\end{remark}

We shall now prove Theorem \ref{MR} (b). Let $g$ be a  positive 
natural number. Let $B=\Spec\ \mC$ and choose 
an elliptic curve $E$ over $\mC$. Let $\CL_E:=\CO(O_E)$ be the 
line bundle associated to the zero-section. Consider 
$E$ as an isotrivial abelian scheme over $H_{1,1}$. 
Consider the abelian scheme 
$\CA:=\CB_{1,1}\times_{H_{1,1}}E^{g-1}$ over $S=H_{1,1}$ and the line 
bundle $\CL:=\CL_{\CB_{1,1}}\boxtimes\CL_E^{\boxtimes(g-1)}$ 
on $\CA$ (here $\boxtimes$ refers to the exterior tensor product). It is a consequence of Proposition \ref{exactor} and Lemma \ref{kun} that 
$\Delta(\CL)$ is of exact order $12$ in $\Pic(S)$.

\begin{acknowledgements} 
The authors wish to thank 
J.-B. Bost for interesting E-mails exchanges related to the 
properties of the bundle $\Delta(\CL)$. They are grateful 
to C.-L. Chai for reading an early version of the text and for providing them with various 
bibliographical and historical references. The authors thanks also go to the referee 
for his very useful suggestions. Finally, they would like 
to thank Ben Moonen for encouraging them to look for the proof 
of the fact that Theorem \ref{MR} (a) is sharp. 
\end{acknowledgements}

\end{document}